\documentclass[11pt,a4paper]{amsart}
\usepackage{amssymb,amsmath}
\usepackage{mathpazo}

\headsep=1cm \numberwithin{equation}{section}
\newtheorem{theorem}{Theorem}[section]
\newtheorem{lemma}[theorem]{Lemma}
\newtheorem{proposition}[theorem]{Proposition}
\newtheorem{corollary}[theorem]{Corollary}

\theoremstyle{definition}
\newtheorem{definition}[theorem]{Definition}
\theoremstyle{remark}
\newtheorem{remark}[theorem]{Remark}
\newtheorem{example}[theorem]{Example}

\newcommand{\Ass}{\operatorname{Ass}}
\newcommand{\im}{\operatorname{im}}
\newcommand{\grade}{\operatorname{grade}}

\newcommand{\Spec}{\operatorname{Spec}}

\newcommand{\ara}{\operatorname{ara}}

\newcommand{\rad}{\operatorname{rad}}
\newcommand{\cd}{\operatorname{cd}}

\newcommand{\Ht}{\operatorname{ht}}
\newcommand{\pd}{\operatorname{pd}}
\newcommand{\Gpd}{\operatorname{Gpd}}
\newcommand{\V}{\operatorname{V}}

\newcommand{\id}{\operatorname{id}}
\newcommand{\Ext}{\operatorname{Ext}}
\newcommand{\Supp}{\operatorname{Supp}}

\newcommand{\Hom}{\operatorname{Hom}}

\newcommand{\Ann}{\operatorname{Ann}}

\newcommand{\depth}{\operatorname{depth}}

\newcommand{\Coass}{\operatorname{Coass}}

\newcommand{\lo}{\longrightarrow}
\newcommand{\fm}{\frak{m}}
\newcommand{\fp}{\frak{p}}

\newcommand{\fa}{\frak{a}}
\newcommand{\fb}{\frak{b}}

\begin{document}

\author[Asgharzadeh, Divaani-Aazar  and Tousi ]{Mohsen Asgharzadeh,
Kamran Divaani-Aazar and Massoud Tousi}

\title[The finiteness dimension of local cohomology ... ]
{The finiteness dimension of local cohomology modules and its dual
notion}

\address{M. Asgharzadeh, Department of Mathematics, Shahid Beheshti
University, Tehran, Iran-and-Institute for Studies in Theoretical
Physics and Mathematics, P.O. Box 19395-5746, Tehran, Iran.}
\email{asgharzadeh@ipm.ir}
\address{K. Divaani-Aazar, Department of Mathematics, Az-Zahra
University, Vanak, Post Code 19834, Tehran, Iran-and-Institute for
Studies in Theoretical Physics and Mathematics, P.O. Box 19395-5746,
Tehran, Iran.} \email{kdivaani@ipm.ir}
\address{M. Tousi, Department of Mathematics, Shahid Beheshti University,
Tehran, Iran-and-Institute for Studies in Theoretical Physics and
Mathematics, P.O. Box 19395-5746, Tehran, Iran.}
\email{mtousi@ipm.ir}


\keywords{Artinianess dimension, finiteness dimension, generalized
local cohomology, local cohomology.\\
The second author was supported by a grant from IPM (No. 86130114).}

\begin{abstract} Let $\fa$ be an ideal of a commutative Noetherian
ring $R$ and $M$ a finitely generated $R$-module. We explore the
behavior of the two notions $f_{\fa}(M)$, the finiteness dimension
of $M$ with respect to $\fa$, and, its dual notion $q_{\fa}(M)$, the
Artinianess dimension of $M$ with respect to $\fa$. When $(R,\fm)$
is local and $r:=f_{\fa}(M)$ is less than $f_{\fa}^{\fm}(M)$, the
$\fm$-finiteness dimension of $M$ relative to $\fa$, we prove that
$H^r_{\fa}(M)$ is not Artinian, and so the filter depth of $\fa$ on
$M$ doesn't exceed  $f_{\fa}(M)$. Also, we show that if $M$ has
finite dimension and $H^i_{\fa}(M)$ is Artinian for all $i>t$, where
$t$ is a given positive integer, then $H^t_{\fa}(M)/\fa
H^t_{\fa}(M)$ is Artinian. It immediately implies that if
$q:=q_{\fa}(M)>0$, then $H^q_{\fa}(M)$ is not finitely generated,
and so $f_{\fa}(M)\leq q_{\fa}(M)$.
\end{abstract}

\maketitle

\section{Introduction}

Throughout this paper,  $R$ is a commutative Noetherian ring with
identity and all modules are assumed to be unitary. Let M be a finitely
generated R-module and $\fa$ an ideal of $R$. The notion
$f_{\fa}(M)$, the finiteness dimension of $M$ relative to $\fa$, is
defined to be the least integer $i$ such that $H^{i}_{\fa}(M)$ is
not finitely generated if there exist such $i$'s and $+\infty$
otherwise. Here $H^{i}_{\fa}(M)$ denotes the $i^{th}$ local
cohomology module of $M$ with respect to $\fa$. As a general
reference for local cohomology, we refer the reader to the text book
\cite{BS}.  Hartshorne \cite{Har} has defined the notion $q_{\fa}(R)$ as the
greatest integer $i$ such that $H^{i}_{\fa}(R)$ is not Artinian.
Dibaei and Yassemi \cite{DY} extended this notion to arbitrary
$R$-modules, to the effect that for any $R$-module $N$ they defined
$q_{\fa}(N)$ as the greatest integer $i$ such that $H^i_{\fa}(N)$ is
not Artinian. Among other things, they showed that if $M$ and $N$
are two finitely generated $R$-modules such that $M$ is supported in
$\Supp_RN$, then $q_{\fa}(M)\leq q_{\fa}(N)$.

Our objective in this paper is to investigate the notions
$f_{\fa}(M)$ and $q_{\fa}(M)$ more closely. Let $M$ and $\fa$ be as above.
By \cite[Theorem 1.2]{AKS}, the $R$-module
$\Hom_R(R/\fa,H^{f_{\fa}(M)}_{\fa}(M))$ is finitely generated.
This easily concludes that $H^{f_{\fa}(M)}_{\fa}(M)$ has finitely many
associated primes, see \cite{BL} and \cite{KS}. In Section 3, we
investigate the dual statements for $H^{q_{\fa}(M)}_{\fa}(M)$. We prove
that $H^{q_{\fa}(M)}_{\fa}(M)/\fa H^{q_{\fa}(M)}_{\fa}(M)$ is Artinian.
However, we give an example to show that the set of coassociated prime
ideals of $H^{q_{\fa}(M)}_{\fa}(M)$ might be infinite, see Example 3.5
below. As  an immediate application, we deduce that if $q:=q_{\fa}(M)>0$,
then $H_{\fa}^{q}(M)$ is not finitely generated. In particular, if
$q_{\fa}(M)>0$, then it follows that $f_{\fa}(M)\leq q_{\fa}(M)$.
This leads one to conjecture that $H^{f_{\fa}(M)}_{\fa}(M)$ is not
Artinian. As can be seen easily, this is not true in general (see
Example 2.6 i) below). But,
in this regard, we prove that if $M$ is a module over a local ring
$(R,\fm)$ such that either $r:=\grade(\fa,M)<\depth M$ or
$r:=f_{\fa}(M)<f_{\fa} ^{\fm}(M)$, then $H^{r}_{\fa}(M)$
is not Artinian. In particular, in both cases we can immediately
conclude that  $f_{\fa}(M)\geq f-\depth(\fa,M)$. (For the definitions
of the notions $f_{\fa}^{\fm}(M)$ and $f-\depth(\fa,M)$,
see the paragraph preceding Theorem 2.5 below.)

To have the most generality, we will present our results for generalized
local cohomology modules, the notion which has been introduced by Herzog
\cite{Her} in 1974. For two $R$-modules $M$ and $N$, the $i^{th}$
generalized local cohomology of $M$ and $N$ with respect to $\fa$ is
defined by  $H^{i}_{\fa}(M,N):=\underset{n}{\varinjlim}\Ext^{i}_{R}
(M/\fa^{n}M, N)$. Using the notion of generalized local cohomology
modules, we can define the finiteness (resp. Artinianess) dimension
of a pair $(M,N)$ of finitely generated $R$-modules relative to the
ideal $\fa$ by
$$f_{\fa}(M,N):=\inf\{i:H^{i}_{\fa}(M,N) \text{ is not finitely
generated}\}$$ (resp. $$q_{\fa}(M,N):=\sup\{i:H^{i}_{\fa}(M,N)
\text{ is not Artinian}\}),$$ with the usual convention that the
infimum (resp. supremum) of the empty set of integers is interpreted
as $+\infty$ (resp. $-\infty$). We study the behavior of $f_{\fa}(M,N)$
and $q_{\fa}(M,N)$  under changing one of the $\fa,M$  and  $N$, when we
fixed the two others. Correspondingly, in Section 2 we prove the following:
\begin{enumerate}
\item[i)] if
$\Supp_R M\subseteq\Supp_R N$, then  $f_{\fa}(M,L)\geq f_{\fa}(N,L)$
for any finitely generated $R$-module $L$, and
\item[ii)] if $\pd N<\infty$, then $f_{\fa}(M,N)\geq f_{\fa}(M,R)
-\pd N$.
\end{enumerate}
Finally, for any two finitely generated $R$-modules $M$ and $N$, we
establish the inequality $$q_{\fa}(M,N)\leq\Gpd_NM+q_{\fa}(M\otimes_RN),$$
where $\Gpd_NM:=\sup\{i:\Ext_R^{i}(M,N)\neq 0\}$.

\section{The finiteness dimension of modules}

The purpose of this section is to examine the behavior of the notion
$f_{\fa}(M)$ more closely. Let's start this section by recording the
following theorem.

\begin{theorem} Let $L,M$ and $N$ be finitely generated $R$-modules and
$\fa$ an ideal of $R$. If $\Supp_R M\subseteq\Supp_R N$, then
$f_{\fa}(M,L)\geq f_{\fa}(N,L)$. In
particular, if $\Supp_R N=\Supp_R M$, then
$f_{\fa}(M,L)=f_{\fa}(N,L)$.
\end{theorem}

{\bf Proof.} It is enough to show that $H_{\fa}^{i}(M,L)$ is a
finitely generated for all $i<f_{\fa}(N,L)$ and all finitely
generated $R$-modules $M$ such that $\Supp_R M\subseteq \Supp_R N$.
To this end, we argue by induction on $i$. We have
$H^0_{a}(M,L)\cong \Hom_R(M,\Gamma_{\fa}(L))$, and so the assertion
is clear for $i=0$. Now, assume that $i>0$ and that the claim has
been proved for $i-1$. By Gruson's theorem (see e.g. [V, Theorem
4.1]), there is a chain
$$0=M_{0}\subseteq M_{1}\subseteq \dots \subseteq M_{\ell}=M$$ of submodules
of $M$ such that each of the factors $M_{j}/M_{j-1}$ is a
homomorphic image of a direct sum of finitely many copies of $N$. In
view of the long exact sequences of generalized local cohomology
modules that induced by the short exact sequences
$$0\longrightarrow M_{j-1}\longrightarrow M_{j}\longrightarrow
M_{j}/M_{j-1}\longrightarrow 0,$$ $j=1,\cdots,\ell$, it suffices to
treat only the case ${\ell}=1$. So, we have an exact sequence
$$0\longrightarrow K\longrightarrow \oplus_{i=1}^tN\longrightarrow
M\longrightarrow 0,$$ where $t\in\mathbb{N}$ and $K$ is a finitely
generated $R$-module. This induces the long exact sequence
$$\cdots\longrightarrow H^{i-1} _{\fa}(K,L) \longrightarrow
H^{i}_{a}(M,L)\longrightarrow
H^{i}_{\fa}(\oplus_{i=1}^tN,L)\longrightarrow \cdots .$$ By the
induction hypothesis, $H^{i-1}_{\fa}(K,L)$ is finitely generated.
Also, $H^{i}_{a}(\oplus_{i=1}^tN,L)\cong
\oplus_{i=1}^tH^{i}_{a}(N,L)$ is finitely generated, because
 $i<f_{\fa}(N,L)$. Hence $H^i_{\fa}(M,L)$ is finitely generated.
$\Box$

\begin{corollary} Let $\fa$ be an ideal of  $R$ and $L,M$ and $N$
finitely generated $R$-modules.
\begin{enumerate}
\item[i)] If  $0\longrightarrow
L\longrightarrow M\longrightarrow N\longrightarrow 0$ is an exact
sequence, then for any finitely generated $R$-module $C$, we have
$f_{\fa}(M,C)=\inf\{ f_{\fa}(L,C),f_{\fa}(N,C)\}.$
\item[ii)] $f_{\fa}(M)=\inf\{ f_{\fa}(C,M):\emph{$C$ is a finitely
generated $R$-module}\}$.
\item[iii)]$f_{\fa}(M,L)=\inf\{ f_{\fa}(R/\fp,L):\fp\in\Supp_R M\}$.
\item[iv)] If $\Supp_R(\frac{M}{\Gamma_{\fa}(M)})\subseteq\Supp_R(\frac{N}
{\Gamma_{\fa}(N)})$, then $f_{\fa}(M,L)\geq f_{\fa}(N,L)$.
\end{enumerate}
\end{corollary}

{\bf Proof.} i) is clear.

ii) Since $H^i_{\fa}(R,M)\cong H^i_{\fa}(M)$, it follows that
$f_{\fa}(R,M)= f_{\fa}(M)$. Now, the claim is clear by Theorem 2.1

iii) Set $K:=\oplus_{\fp \in \Ass_RM}R/\fp$. Then $K$ is finitely
generated and $\Supp_RK=\Supp_RM$. So, by i) and Theorem 2.1, we
deduce that
$$\begin{array}{ll} f_{\fa}(M,L)&=f_{\fa}(K,L)\\
&=\inf\{f_{\fa}(R/\fp,L):\fp\in \Ass_RM\}\\
&=\inf\{f_{\fa}(R/\fp,L):\fp\in \Supp_RM\}.
\end {array}$$

iv) \cite[Lemma 2.11]{DH} implies that $H^i_{\fa}(C,L)\cong
\Ext_R^i(C,L)$ for all $\fa$-torsion $R$-modules $C$ and all $i$.
So, the exact sequence $0\lo \Gamma_{\fa}(M)\lo M \lo
\frac{M}{\Gamma_{\fa}(M)}\lo 0$ implies the long exact sequence
$$\cdots \lo \Ext_R^{i-1}(\Gamma_{\fa}(M),L)\lo H^i_{\fa}
(\frac{M}{\Gamma_{\fa}(M)},L)\lo H^i_{\fa}(M,L)\lo
\Ext_R^i(\Gamma_{\fa}(M),L)\lo \cdots .$$ This yields that
$f_{\fa}(M,L)=f_{\fa}(\frac{M}{\Gamma_{\fa}(M)},L)$ and similarly we
have $f_{\fa}(N,L)=f_{\fa}(\frac{N}{\Gamma_{\fa}(N)},L)$. Now, the
claim becomes clear by Theorem 2.1.  $\Box$

\begin{example} Let $L$, $M$ and $N$ be finitely generated $R$-modules
such that $\Supp_RM=\Supp_RN$. In Theorem 2.1, we saw that $f_{\fa}(M,L)
=f_{\fa}(N,L)$ for any ideal $\fa$ of $R$. One may ask whether the equality
$f_{\fa}(L,M)=f_{\fa}(L,N)$ holds too. This would not be the case. To
see this, let $(R,\fm)$ be a local ring with $\depth R>1$ and $N'$
a 1-dimensional Cohen-Macaulay $R$-module. Set $L=M:=R$ and $N:=N'\oplus R$.
Then $\Supp_RM=\Supp_RN$. Now, by \cite [Remark 2.5]{Hel2},
$H^1_{\fm}(N)\cong H^1_{\fm}(N')$ is not finitely generated. So
$f_{\fm}(L,N)=f_{\fm}(N)=1$, while $f_{\fm}(L,M)=f_{\fm}(R)>1$.
\end{example}

If we fix the ideal $\fa$ and the $R$-module $M$, then we can't say
so much about $f_{\fa}(M,\cdot)$. However, we have the following
result.

\begin{proposition} Let $\fa$ be an ideal of  $R$ and  $M, N$ two
finitely generated $R$-modules such that $\pd N<\infty$. Then
$f_{\fa}(M,N)\geq f_{\fa}(M,R)-\pd N$.
\end{proposition}

{\bf Proof.} We use induction on $n:=\pd N$. If $n=0$, then there is
a nothing to prove. Now, assume that $n>0$ and that the assertion
holds for $n-1$. We can construct an exact sequence
$$0\longrightarrow L\longrightarrow F\longrightarrow
N\longrightarrow 0$$  of finitely generated $R$-modules such that
$F$ is free and $\pd L=n-1$. By the induction hypothesis,
$f_{\fa}(M,L)\geq f_{\fa}(M,R)-n+1$. Let $i<f_{\fa}(M,R)-n$. Then,
it follows from the exact sequence
$$H_{\fa}^i(M,F)\longrightarrow H_{\fa}^i(M,N)\longrightarrow
H_{\fa}^{i+1}(M,L)$$ that $H_{\fa}^i(M,N)$ is finitely generated.
Hence $f_{\fa}(M,N)\geq f_{\fa}(M,R)-n$, as required. $\Box$

Let $\fa\subseteq \fb$ be  ideals of $R$. Recall that
\cite[Definition 5.3.6]{S} for a not necessary finitely generated
$R$-module $M$, $\grade(\fa,M)$ is defined by
$$\grade(\fa,M):=\inf\{i\in \mathbb{N}_0:\Ext_R^i(R/\fa,M)\neq
0\}.$$ Also, recall that $f_{\fa}^{\fb}(M)$, the $\fb$-finiteness
dimension $M$ relative to $\fa$, is defined by
$$f_{\fa}^{\fb}(M):=\inf\{i\in \mathbb{N}_0:\fb \nsubseteq
\rad(\Ann_RH^i_{\fa}(M))\}.$$ If $M$ is finitely generated, then
\cite[Proposition 9.1.2]{BS} implies that $f_{\fa}^{\fa}(M)=
f_{\fa}(M)$. Also, we remind the reader that for a finitely
generated $R$-module $M$ over a local ring $(R,\fm)$, the filter
depth of $\fa$ on $M$ is defined as the length of any maximal
$M$-filter regular sequence in $\fa$ and denoted by
$f-\depth(\fa,M)$, see \cite[Defintion 3.3]{LT}. By \cite[Theorem
3.1]{M1}, it is known that $f-\depth(\fa,M)= \inf\{i\in
\mathbb{N}_0:H^i_{\fa}(M) \text{ is not Artinian}\}.$

\begin{theorem} Let $\fa$ be an ideal of the local ring $(R,\fm)$
and $M$ an $R$-module. Then, the  following holds:
\begin{enumerate}
\item[i)] Assume that $r:=f_{\fa}^{\fa}(M)<f_{\fm}^{\fa}(M)$. Then $H^{r}
_{\fa}(M)$ is not Artinian. Moreover, if $M$ is finitely generated,
then $f-\depth(\fa,M)\leq f_{\fa}(M)$.
\item[ii)] Assume that $r:=\grade(\fa,M)<\depth M$. Then $H^{r}
_{\fa}(M)$ is not Artinian. Moreover, if $M$ is finitely generated,
then $\grade(\fa,M)=f-\depth(\fa,M)$.
\end{enumerate}
\end{theorem}

{\bf Proof.} i) We argue by induction on $r$. Let $r=0$. If $H^{0}
_{\fa}(M)$ is Artinian, then $$H^{0} _{\fa}(M)\cong
H^{0}_{\fm}(H^{0} _{\fa}(M) )\cong H^{0}_{\fm}(M).
$$ Since $0<f_{\fm}^{\fa}(M)$, there is an integer $n$ such that
$\fa^{n}H^0_{\fa}(M)\cong \fa^{n}H^0_{\fm}(M)=0$. So
$0<f_{\fa}^{\fa}(M)$ and we achieved at a contradiction.

Now, assume that $r>0$. Then there exists an integer
$n\in\mathbb{N}$ such that $\fa^{n}H^0_{\fa}(M)=0$. Hence the
argument \cite[Remark 1.3 ii)]{BRS} yields that
$f_{\fm}^{\fa}(M)=f_{\fm}^{\fa}(M/\Gamma _{\fa}(M))$ and
$f_{\fa}^{\fa}(M)=f_{\fa}^{\fa}(M/\Gamma_{\fa}(M))$. Thus without
loss of generality, we may and do assume that $H^0_{\fa}(M)=0$. Now,
we apply Melkersson's technic \cite{M2}, so let $E$ be an injective
envelope of $M$ and $N:=E/M$. Then
$\Gamma_{\fa}(E)=\Gamma_{\fm}(E)=0$. From the exact sequence of
local cohomology modules induced by
$$0\longrightarrow M\longrightarrow E\longrightarrow N
\longrightarrow0,$$ we obtain that $H^{i} _{\fa}(N)\cong H^{i+1}
_{\fa}(M)$ and $H^{i} _{\fm}(N)\cong H^{i+1} _{\fm}(M)$ for all
$i\geq0$. Hence $f_{\fm}^{\fa}(N)=f_{\fm}^{\fa}(M)-1$ and
$f_{\fa}^{\fa}(N) =f_{\fa}^{\fa}(M)-1$ and the claim follows by the
induction hypothesis.

ii) The proof is similar to the proof of i). Note that by
\cite[Proposition 5.3.15]{S}, $\grade(\fa,M)=\inf\{i\in
\mathbb{N}_0:H^i_{\fa}(M)\neq 0\}.\ \   \Box$

\begin{example}
\begin{enumerate}
\item[i)] The assumption $f_{\fa}^{\fa}(M)<f_{\fm}^{\fa}(M)$ is crucial
in Theorem 2.5 i). To realize this, let $M$ be a Cohen-Macaulay
$R$-module of positive dimension and $\fa=\fm$. Then $f_{\fa}^{\fa}(M)=\dim M$ and $H^{\dim M}_{\fa}(M)$ is Artinian.
\item[ii)] Also, the assumption $\grade(\fa,M)<\depth M$ cannot be dropped in
Theorem 2.5 ii). In fact, if $M$ is a finitely generated $R$-module and $\fa=
\fm$, then $\grade(\fa,M)=\depth M$ and all local cohomology modules $H^{i}_{\fa}(M)$ are Artinian.
\end{enumerate}
\end{example}

\begin{remark}
\begin{enumerate}
\item[i)] Let $M$ be a finitely generated $R$-module such that
$\Ht_{M}\fa>0$, then we have the inequality $f_{\fa}(M)\leq
\Ht_{M}\fa$. To see this, first note that we may  assume that $\fa
M\neq M$. Let $\fp\in \Supp_R (M/\fa M)$ be such that
$\Ht_{M}\fa=\Ht_{M}\fp$, and set $t=\Ht_{M}\fa$. Then, because of
the natural isomorphisms
$$H^t_{\fa}(M)_{\fp}\cong H^t_{(\fa+\Ann_RM)R_{\fp}}(M_{\fp})
\cong H^t_{\fp R_{\fp}}(M_{\fp})$$ and \cite[Remark 2.5]{Hel2}, it
follows that $H^{t}_{\fa}(M)_{\fp}$ is not a finitely generated
$R_{\fp}$-module.  Hence the $R$-module $H^t_{\fa}(M)$ is not
finitely generated and consequently, $f_{\fa}(M)\leq \Ht_{M}\fa$.
\item[ii)] From the definition of $f_{\fa}(M)$, it becomes clear
that $f_{\fa}(M)\geq \grade(\fa,M)$. There are some cases in
which the equality holds. For example, let $M$ be a  Cohen Macaulay
$R$-module such that $\Ht_{M}\fa>0$, then by i), $f_{\fa}(M)\leq
\Ht_{M}\fa= \grade(\fa,M)$, and so $f_{\fa}(M)=\grade(\fa,M)$.
\item[iii)] If $(R,\fm)$ is a regular local ring with $\dim R>1$, then for any
integer $0<n<\dim R$, there exists an $n$-dimensional finitely
generated $R$-module $M$ such that $\depth M=0$ and $f_{\fm}(M)=\dim
M$. This holds, because by \cite[Ex. 6.2.13]{BS}, there is a
finitely generated $R$-module $M$ such that $H^i_{\fm}(M)\neq 0$ if
and only if $i$ is either $0$ or $n$. So,  $\depth M=0$ and by
\cite[Exersice 6.1.6]{BS}, $f_\fm(M)=n=\dim M$.
\item[iv)] Let $(R,\fm)$ be a local ring, $\fa$ a proper ideal of $R$ and
$M,N$ finitely generated $R$-modules. Set
$\fa_M:=\Ann_R(\frac{M}{\fa M})$. By \cite[Proposition 5.5]{B}, it
follows that
$$\inf\{i\in \mathbb{N}_0:H^i_{\fa}(M,N)\neq 0\}=\inf\{i\in
\mathbb{N}_0:H^i_{\fa_M}(N)\neq 0\}(=\grade(\fa_M,N)).$$
Also, \cite[Theorem 2.2]{CT} yields that $$\inf\{
i\in\mathbb{N}_0:H^i_{\fa}(M,N) \text{ is not Artinian}\}=\inf\{
i\in\mathbb{N}_0:H^i_{\fa_M}(N) \text{ is not Artinian}\}$$$(=
f-\depth(\fa_M,N)).$ Having these facts in mind, one might ask
whether $f_{\fa}(M,N)= f_{\fa_M}(N)$. This is not necessarily true.
For instance, let $M:=R/\fm$, $\fa:=\fm$ and $N$ be any non-Artinian
finitely generated $R$-module. Then
$$H^{i}_{\fa}(M,N)=\underset{n}{\varinjlim}\Ext^{i}_{R}(M/\fa^{n}M,N)
\cong \Ext^{i}_{R}(R/\fm, N)$$ is finitely generated for all $i$,
and so $f_{\fa}(M,N)=+\infty$, while $f_{\fa_M}(N)\leq \dim
N<\infty$.
\end{enumerate}
\end{remark}

\section{The Artinian dimension of modules}

In this section, we focus on the invariant $q_{\fa}(M,N)$. In
\cite[Definition 2.1]{DM}, the authors call an $R$-module $N$ weakly
Laskerian if any quotient of $N$ has finitely many associated prime
ideals. In the sequel, for a not necessarily finitely generated $R$-module
$N$, by dimension of $N$, we mean the dimension of $\Supp_RN$.

\begin{theorem} Let $\fa$ be an ideal of  $R$,  $M$ a finitely
generated $R$-module of finite projective  dimension and $N$ a
weakly Laskerian $R$-module of finite dimension. Let $t>\pd M$ be an
integer such that $H_{\fa}^{j}(M,N)$ is Artinian for all $j>t$. Then
$H_{\fa}^{t}(M,N)/\fa H_{\fa}^{t}(M,N)$ is Artinian.
\end{theorem}

{\bf Proof.}  We use induction on $n:=\dim N$. By \cite[Theorem
2.5]{DH}, it follows  that $H_{\fa}^i(M,L)=0$ for all finitely generated
$R$-modules $L$ and all $i>\pd M+\dim L$. So, since the functor
$H_{\fa}^i(M,\cdot)$ commutes with direct limits, it turns out that
$H_{\fa}^i(M,N)=0$ for all $i>\pd M+\dim N$. Thus the claim clearly holds
for $n=0$.

Now, assume that $n>0$ and that the claim holds for $n-1$. Since
$t>\pd M$, in view of the long exact sequence of generalized local
cohomology modules that is induced by the exact sequence $$0\lo
\Gamma_{\fa}(N)\lo N \lo N/\Gamma_{\fa}(N)\lo 0,$$ we may assume
that $N$ is $\fa$-torsion free. Note that since the functor
$H_{\fa}^i(M,\cdot)$ is the $i^{th}$ right derived functor of the functor
$\Hom_R(M,\Gamma_{\fa}(\cdot))$ and $\Gamma_{\fa}(N)$ possesses an
injective resolution consisting of $\fa$-torsin injective $R$-modules,
it follows that $H_{\fa}^i(M, \Gamma_{\fa}(N))\cong
\Ext_R^i(M,\Gamma_{\fa}(N))$ for all $i$. Take $x\in\fa\setminus \bigcup _{\fp\in \Ass_RN}\fp$. Then
$N/xN$ is weakly Laskerian and $\dim N/xN\leq n-1$. The exact sequence
$$0\longrightarrow N \stackrel{x}\longrightarrow N\longrightarrow
N/xN\longrightarrow 0$$ implies the following long exact sequence of
generalized local cohomology modules
$$\cdots \lo H_{\fa}^j(M,N)\stackrel{x}\lo H_{\fa}^j(M,N)\lo
H_{\fa}^j(M,N/xN)\lo \cdots .$$ It yields that $H^{j}_{\fa}(M,N/xN)$
is Artinian for all $j>t$. Thus $\frac{H^{t}_{\fa}(M,N/xN)}{\fa
H^{t}_{\fa}(M,N/xN)}$ is Artinian by the induction hypothesis.

Now, consider the exact sequence $$
H^{t}_{\fa}(M,N)\stackrel{x}\longrightarrow
H^{t}_{\fa}(M,N)\stackrel{f}\longrightarrow
H^{t}_{\fa}(M,N/xN)\stackrel{g}\longrightarrow H^{t+1}_{\fa}(M,N),$$
which induces the following two exact sequences
$$0\longrightarrow \im f \longrightarrow H^{t}_{\fa}(M,N/xN)
\longrightarrow\im g \longrightarrow 0$$ and $$
H^{t}_{\fa}(M,N)\stackrel{x}\longrightarrow
H^{t}_{\fa}(M,N)\longrightarrow \im f\longrightarrow 0.$$ Therefore,
we can obtain the following two exact sequences:
$$Tor^{R}_{1}(R/\fa,\im g )\lo \im f/\fa \im f\lo
\frac{H^{t}_{\fa}(M,N/xN)}{\fa H^{t}_{\fa}(M,N/xN)}\lo \im g/\fa \im
g\lo 0   (*)$$ and
$$\frac{H^{t}_{\fa}(M,N)}{\fa H^{t}_{\fa}(M,N)}\stackrel{x}
\longrightarrow \frac{H^{t}_{\fa}(M,N)}{\fa
H^{t}_{\fa}(M,N)}\longrightarrow \im f/\fa \im f\longrightarrow0.$$
Since $x\in \fa$, from later exact sequence, we deduce that
$\frac{H^{t}_{\fa}(M,N)}{\fa H^{t}_{\fa}(M,N)}\cong \im f/\fa \im
f$. Now, since $Tor^{R}_{1}(R/\fa,\im g)$ and
$\frac{H^t_{\fa}(M,N/xN)}{\fa H^t_{\fa}(M,N/xN)}$ are Artinian, the
claim follows by $(*)$. $\Box$

The following Corollary improves \cite[Theorem 4.7]{DN}.

\begin{corollary} Let $\fa$ be an ideal of $R$
and $N$ a weakly Laskerian $R$-module of finite dimension. Let $t$
be a positive integer. If $H_{\fa}^i(N)$ is Artinian for all $i>t$,
then $H_{\fa}^t(N)/\fa H_{\fa}^t(N)$ is Artinian.
\end{corollary}

\begin{corollary} Let $\fa$ be an ideal of $R$,
$M$ a finitely generated  $R$-module of finite projective dimension
and $N$ a weakly Laskerian $R$-module of finite dimension.  If
$q:=q_{\fa}(M,N)>\pd M$, then  $H_{\fa}^{q}(M,N)$ is not finitely
generated, and so $f_{\fa}(M,N)\leq q_{\fa}(M,N) $. In particular,
if If $r:=q _{\fa}(N)>0$, then $H_{\fa}^{r}(N)$ is not finitely
generated and so $f_{\fa}(N)\leq q_{\fa}(N)$.
\end{corollary}

{\bf Proof.} Contrary, assume that $H_{\fa}^{q}(M,N)$ is finitely
generated. Then, there exists an integer $n$ such that $\fa^{n}
H_{\fa}^q(M,N)=0$. Since $\fa^{n}$ and $\fa$ have the same radical,
it turns out that $H_{\fa ^{n}}^i(M,N)\cong H_{\fa}^i(M,N)$ for all
$i$. Thus, Theorem 3.1 yields that $H_{\fa}^{q}(M,N)\cong
\frac{H^q_{\fa ^{n}}(M,N)}{\fa^{n} H^q_{\fa ^{n}}(M,N)}$ is
Artinian, and so we achieved at a contradiction. $\Box$

In the sequel, we use the notion of cohomological
dimension. Recall that for an $R$-module $N$, the cohomological
dimension of $N$ with respect to $\fa$ is defined by
$\cd(\fa,N):=\sup\{i\in \mathbb{N}_0:H_{\fa}^i(N)\neq 0\}$. Also,
recall that the arithmetic rank $\ara(\fa)$ of the ideal $\fa$ is
the least number of elements of $R$ required to generate an ideal
which has the same radical as $\fa$. By \cite[Corollary 3.3.3 and
Theorem 6.1.2]{BS}, it turns out
that $\cd(\fa,N)\leq \min\{\ara(\fa),\dim N\}$.

\begin{example}
\begin{enumerate}
\item[i)] In  Corollary 3.2, the positivity assumption on $t$ is
really necessary. To see this, let $(R,\fm)$ be a local ring and
consider the weakly Laskerian $R$-module
$M:=\bigoplus_{\mathbb{N}}R/\fm$. We have $H_{\fm}^i(M)=0$ for all
$i>0$, but $\frac{H_{\fm}^{0}(M)}{\fm
H_{\fm}^0(M)}=\bigoplus_{\mathbb{N}}R/\fm$ is not Artinian.
\item[ii)] In Corollary 3.2, if $t:=\cd(\fa,N)$, then it can be
seen easily that $\frac{H_{\fa}^{t}(N)}{\fa H_{\fa}^{t}(N)}\cong
H_{\fa}^{t}(N/\fa N)=0$. But, in general $\frac{H_{\fa}^{t}(N)}{\fa
H_{\fa}^{t}(N)}$ even might not be finitely generated. To see this,
let $R:=k[[X_1,X_2,X_3,X_4]]$, where $k$ is a field. Let $\fp_{1}:=(X_{1},X_{2}),
 \fp_{2}:=(X_{3},X_{4})$ and $\fa:=\fp_{1}\cap \fp_{2}$. By
\cite[Example 3]{Har}, one has that $q_{\fa}(R)=2$ and
$H_{\fa}^{2}(R)\cong H_{\fp_{1}}^{2}(R)\oplus H_{\fp_{2}}^{2}(R)$.
Now consider the following isomorphisms
$$\begin{array}{ll}
\frac{H_{\fa}^{2}(R)}{\fa H_{\fa}^{2}(R)}&\cong
\frac{H_{\fp_{1}}^{2}(R)}{\fa H_{\fp_{1}}^{2}(R)}\oplus
\frac{H_{\fp_{2}}^{2}(R)}{\fa H_{\fp_{2}}^{2}(R)} \\&\cong
H_{\fp_{1}}^{2}(R/\fa)\oplus H_{\fp_{2}}^{2}(R/\fa).
\\
\end{array}
$$
By the Hartshorne-Lichtenbaum Vanishing Theorem,
$H_{\fp_{1}}^2(R/\fa)\neq 0$. Therefore $\cd(\fp_{1},R/\fa)=2$, and
so by \cite[Remark 2.5]{Hel2}, $H_{\fp_{1}}^{2}(R/\fa)$ is not
finitely generated. Consequently, $\frac{H_{\fa}^{2}(R)}{\fa
H_{\fa}^{2}(R)}$ is not finitely generated.
\end{enumerate}
\end{example}

Recall that for an $R$-module $M$, $\Coass_RM$, the set of
coassociated prime ideals of $M$, is defined to be the set of all
prime ideals $\fp$ of $R$ such that $\fp=\Ann_RL$ for some Artinian
quotient $L$ of $M$. In the case $(R,\fm)$ is local, it is known
that $\Coass_RM=\Ass_R(\Hom_R(M,E(R/\fm)))$ (see e.g.  \cite[Theorem
1.7]{Y}).

\begin{example}
Let $\fa$ be an ideal of $R$ and $M$ a finitely generated
$R$-module. Let $t$ be a natural integer. Assume that $H^i_{\fa}(M)$
is finitely generated for all $i<t$. Then by \cite[Theorem 1.2]{AKS},
$\Hom_R(R/\fa,H^t_{\fa}(M))$ is finitely generated. (This can be
viewed as dual of Corollary 3.2.) Since $H^t_{\fa}(M)$ is
$\fa$-torsion, it follows that
$\Ass_R(H^t_{\fa}(M))=\Ass_R(\Hom_R(R/\fa,H^t_{\fa}(M)))$, and so
$\Ass_R(H^t_{\fa}(M))$ is finite. Now, assume that $H^i_{\fa}(M)$ is
Artinian for all $i>t$. Then, it is rather natural to expect that the set
$\Coass_R (H^t_{\fa}(M))$ is finite. But, this is not the case.
To see this, let $(R,\fm)$ be an equicharacteristic local ring with
$\dim R>2$. Let $\fp$ be a prime ideal of $R$ of height 2 and take
$x\in \fm-\fp$. Then by \cite[Corollary 2.2.2]{Hel1},
$\Coass_R(H_{(x)}^1(R))=\Spec R\setminus \V((x))$. Since $\Ht \fp=2$,
there are infinitely many prime ideals of $R$ which are contained in $\fp$,
and so $\Coass_R(H_{(x)}^1(R))$ is infinite.
\end{example}

We apply the following lemma in the proof of next theorem.

\begin{lemma} Let $A$ be an Artinian $R$-module and $S$ a
multiplicatively closed subset of $R$. Then as an $R$-module,
$S^{-1}A$ is isomorphic to a submodule of $A$, and so $S^{-1}A$ is
Artinian both as an $R$-module and as an $S^{-1}R$-module.
\end{lemma}

{\bf Proof.} Since $A$ is Artinian, it is supported in finitely many
maximal ideals $\fm_1,\dots, \fm_t$ of $R$. It is known that there is a
natural isomorphism $A\cong \oplus_{i=1}^t\Gamma_{\fm_i}(A)$. We can
assume that there is an integer $1\leq {\ell}\leq t$ such that $S$
doesn't intersect $\fm_1,\dots, \fm_{\ell}$, while $S\cap
\fm_i\neq \emptyset$ for the remaining $i$'s.  For any maximal ideal
$\fm$, one can check that $S^{-1}(\Gamma_{\fm}(A))=0$ if $S\cap
\fm\neq \emptyset$ and that the natural map $\Gamma_{\fm}(A)\lo
S^{-1}(\Gamma_{\fm}(A))$ is an isomorphism otherwise.  Hence, we
have a natural $R$-isomorphism $S^{-1}A\cong
\oplus_{i=1}^{\ell}\Gamma_{\fm_i}(A)$. This shows that as an
$R$-module, $S^{-1}A$ is isomorphic to a submodule of $A$. In
particular, $S^{-1}A$ is an Artinian $R$-module, and so it is also
Artinian as an $S^{-1}R$-module. $\Box$

Part iv) of the next result has been proved by Chu and Tang
\cite[Theorem 2.6]{CT} under the extra assumption that $R$ is local.

\begin{theorem} Let $L,M$ and $N$ be finitely generated $R$-modules and
$\fa\subseteq\fb$ ideals of $R$. Then the following holds:
\begin{enumerate}
\item[i)]$q_{\fb}(M,N)\leq q_{\fa}(M,N)+\ara(\fb/\fa)$.
\item[ii)] If $(R,\fm)$ is local and $\ara(\fm/\fa)\leq 1$, then
$q_{\fa}(M,N)=\sup\{i: \Supp_R(H^i_{\fa}(M,N))\nsubseteqq
\{\fm\}\}$.
\item[iii)] If $S$ is a multiplicatively closed subset of $R$,
then $q_{S^{-1}\fa}(S^{-1}M,S^{-1}N)\leq q_{\fa}(M,N)$.
\item[ iv)] If $\pd M<\infty$ and $\Supp_R N\subseteq\Supp_R L$, then
$q_{\fa}(M,N)\leq q_{\fa}(M,L)$. In particular, if $\Supp_R
N=\Supp_R L$, then $q_{\fa}(M,N)=q_{\fa}(M,L)$.
\item[v)] If $\pd M<\infty$, then $q_{\fa}(M,N)=
\sup\{q_{\fa}(M,R/\fp):\fp\in \Supp_R N\}$.
\end{enumerate}
\end{theorem}

{\bf Proof.} i) Since the generalized local cohomology functors with
respect to ideals with the same radicals are equivalent, we may
assume that there are $x_1,x_2,\dots x_n \in R$ such that
$\fb=\fa+(x_1,x_2,\dots x_n)$. By induction on $n$, it is enough to
show the claim for the case $n=1$. So, let $\fb=\fa+(x)$ for some
$x\in R$. By \cite[Lemma 3.1]{DH}, there is the following long exact
sequence of generalized local cohomology modules
$$\cdots \lo H_{\fa R_x}^{i-1}
(M_x,N_x) \lo H_{\fb}^i(M,N)\lo H_{\fa}^i(M,N)\lo \cdots .$$ Let
$i>q_{\fa}(M,N)+1$. Then $H_{\fa}^i(M,N)$ and $H_{\fa}^{i-1}(M,N)$
are Artinian. Also, by Lemma 3.6, it turns out that $H_{\fa
R_x}^{i-1}(M_x,N_x)\cong H_{\fa}^{i-1}(M,N)_x$ is Artinian as an
$R$-module. Thus $H_{\fb}^i(M,N)$ is Artinian, and so
$q_{\fb}(M,N)\leq q_{\fa}(M,N)+1$.

ii) First of all note that $H_{\fm}^i(M,N)$ is Artinian for all $i$.
Similar to i), we may assume that $\fm=\fa+(x)$ for some $x\in R$.
Clearly, $$t:=\sup\{i: \Supp_R(H^i_{\fa}(M,N))\nsubseteqq
\{\fm\}\}\leq q_{\fa}(M,N).$$ Let $j>t$ be an integer. Then
$H_{\fa}^j(M,N)$ is $\fm$-torsion, and so
$$H_{\fa R_x}^j(M_x,N_x)\cong H_{\fa}^j(M,N)_x=0.$$ Hence, from the
exact sequence
$$\cdots\lo H_{\fm}^j(M,N)\lo H_{\fa}^j(M,N)\lo
H_{\fa R_x}^j(M_x,N_x)\lo \cdots,$$ we conclude that
$H_{\fa}^j(M,N)$ is Artinian. This implies that $q_{\fa}(M,N)\leq
t$.

iii) is clear by Lemma 3.6.

iv) For any finitely generate $R$-module $C$, one has $H_{\fa}^i(M,C)=0$
for all $i>\pd M+\ara(\fa)$, see e.g. \cite [Theorem 2.5]{DH}.
Hence by using decreasing induction on $q_{\fa}(M,L) \leq i\leq
\pd M+\ara(\fa)+1$, one can prove the claim by slight modification
of the proof of Theorem 2.1.

v) this can be deduced by an argument similar to that used in the
proof of Corollary 2.2 iii). $\Box$

\begin{remark} It is known that if $\fa$ is an ideal of a regular
local ring $(R,\fm)$, then (under some mild assumptions on $R$), we have $q_{\fa}(R)=\sup\{i\in
\mathbb{N}_0:\Supp_R(H_{\fa}^i(R))\nsubseteq \{\fm\}\}$. See
\cite[Corollary 2.4]{HS} for the case that the characteristic
of $R$ is prime and \cite[Theorem 2.7]{O} for the
characteristic 0 case. Part ii)
of the above theorem might be considered as a generalization of the
Hartshorne-Speiser and Ogus's results. But, the reader should be
aware that the assumption $\ara(\fm /\fa)\leq 1$ is necessary. To
see this, let $R = k[[U, V, X, Y]]/(UX + VY)$; $k$ a field,
$\fm:=(U,V,X,Y)R$ and $\fa=(U,V)R$. Then
$\Supp_R(H_{\fa}^{i}(R))\subseteq \{\fm\}$
 for all $i\geq 2$ and
$H_{\fa}^{2}(R)$ is not Artinian, see \cite[Theorem 1.1]{HelS} . Note
that $\ara(\fm/\fa)=2$ and $(R,\fm)$ is a complete intersection ring
which is not regular.
\end{remark}

In the sequel, we need the following definition from \cite{DH}.

\begin{definition} Let $M$ and $N$ be finitely generated
$R$-modules. We define {\it projective dimension of $M$ relative to
$N$} by $$\pd_NM:=\sup\{\pd_{R_{\fp}}M_{\fp}:\fp\in \Supp_RM\cap
\Supp_RN \}.$$ Also, we define {\it Gorenstein projective dimension
of $M$ relative to $N$} by
$$\Gpd_NM:=\sup\{i\in \mathbb{N}_0:\Ext_R^i(M,N)\neq 0 \}.$$
\end{definition}

\begin{lemma} Let $M$ and $N$ be finitely generated
$R$-modules.
\begin{enumerate}
\item[i)] $\grade(\Ann _R M,N)\leq \Gpd_NM$.
\item[ii)] If $\pd_NM$ is finite, then
$\Gpd_NM=\pd_NM$. In particular, if $R$ is local and $\pd M$ is
finite, then $\Gpd_NM=\pd M$.
\item[iii)] $\Gpd_NM\leq \min\{\pd M, \id N \}$. In particular, if
either $\pd M$ or $\id N$ is finite, then $\Gpd_NM$ is finite.
\item[iv)] If $\pd_NM$ is finite, then $\bigcup_{i\in \mathbb{N}_0}\Supp_R
(\Ext_R^i(M,N))=\Supp_RM\cap \Supp_RN$, and so  $\dim
(M\otimes_RN)=\sup\{\dim (\Ext_R^i(M,N)):i\in \mathbb{N}_0\}$.
\item[v)] If $\pd_NM$ is finite,  then $q_{\fa}(L,M\otimes_RN)=
\max\{q_{\fa}(L,\Ext_R^i(M,N)):i\in \mathbb{N}_0 \}$ for all
finitely generated $R$-modules $L$ of finite projective dimension.
\item[vi)] If $R$ is local and $\id N<\infty$, then
$\Gpd_NM=\depth R-\depth M$, and if in addition $M$ is maximal
Cohen-Macaulay, then $\Gpd_NM=0$.
\end{enumerate}
\end{lemma}

{\bf Proof.} We only prove v). The other parts are proved in
\cite[Lemma 2.2]{DH}.

By iv) and Theorem 3.7 v), we have
$$
\begin{array}{ll} q_{\fa}(L,M\otimes_RN)
&=\sup\{q_{\fa}(L,R/\fp):\fp\in \Supp_R(M\otimes_RN)\}
\\&=\sup \{q_{\fa}(L,R/\fp):\fp\in \bigcup_{i\in \mathbb{N}_0}
\Supp_R(\Ext_R^i(M,N))\}
\\&=\underset{i\in \mathbb{N}_0}\max(\sup\{q_{\fa}(L,R/\fp):
\fp\in \Supp_R(\Ext_R^i(M,N))\})\\
&=\max\{q_{\fa}(L,\Ext_R^i(M,N)):i\in \mathbb{N}_0 \}. \  \  \Box
\end{array}$$

\begin{theorem} Let $M$ and $N$ be finitely generated $R$-modules.
For any ideal $\fa$ of $R$, the inequality $q_{\fa}(M,N)\leq
\Gpd_NM+q_{\fa}(M\otimes_{R}N)$ holds.
\end{theorem}

{\bf Proof.} We may assume that $\Gpd_NM<\infty$. Let $F(\cdot):=
\Gamma_{\fa}(\cdot)$ and $G(\cdot):=\Hom_R(M,\cdot)$. It is straightforward
to see that $H_{\fa}^i(\Hom_R(M,E))=0$ for any injective $R$-module $E$ and
all $i\geq 1$. So, since $(GF)(\cdot)=\Hom_R(M,\Gamma_{\fa}(\cdot))$, by
\cite [Theorem 11.38]{R}, one has the following Grothendieck's spectral sequence
$$H^i_{\fa}(\Ext_R^j(M,N))\underset{i}\Longrightarrow
H^{i+j}_{\fa}(M,N).$$  Hence for each $n\in\mathbb{N}_{0}$, there
exists a chain
$$(\ast)    \  \ 0=H^{-1}\subseteq H^{0}\subseteq\cdots \subseteq
H^{n}:=H^n_{\fa}(M,N)$$ of submodules of $H^n_{\fa}(M,N)$ such
that $H^i/H^{i-1}\cong E^{i,n-i}_{\infty}$ for all $i=0,\cdots,n$.
Since, $\Ext_R^j(M,N)$ is supported in $\Supp_RM \cap \Supp_RN$, it
follows by \cite[Theorem 3.2]{DY} that
$H^i_{\fa}(\Ext_{R}^{j}(M,N))$ is Artinian for all
$i>q_{\fa}(M\otimes_{R}N)$. If
$n>\Gpd_NM+q_{\fa}(M\otimes_{R}N)$, then either
$i>q_{\fa}(M\otimes_{R}N)$ or  $n-i>\Gpd_NM$.  In each case, it turns out
that $E_2^{i,n-i}$ is Artinian. Therefore, by splitting the
chain $(\ast)$ into short exact sequences, we deduce that
$H^n_{\fa}(M,N)$ is Artinian for all
$n>\Gpd_NM+q_{\fa}(M\otimes_{R}N)$. Note that $ E_{\infty}^{i,n-i}$
is a subquotient of  $E_2^{i,n-i}$ for all $i$. $\Box$

In the sequel, we consider some cases in which the equality holds in
Theorem 3.11.

\begin{corollary} Let $\fa$ be an ideal of $R$ and $M,N$ two finitely
generated $R$-modules.
\begin{enumerate}
\item[i)]If $p:=\Gpd_NM=\grade(\Ann _R M,N)$, then $H^i_{\fa}
(M,N)=H^{i-p}_{\fa}(\Ext^p_R(M,N))$ for all $i$, and consequently
$q_\fa(M,N)=q_\fa (\Ext^p_R(M,N))+p$. Also, $f_\fa(M,N)=f_\fa
(\Ext^p_R(M,N))+p$.
\item[ii)] If $p:=\pd_NM=\grade(\Ann _R M,N)$, then $q_\fa(M,N)
=q_\fa (M\otimes_R N)+p$.
\item[iii)] Let $\fa$ be an ideal of the local ring $(R,\fm)$ and $M,N$ two
finitely generated $R$-modules such that $M$ is  faithful and
maximal Cohen Macaulay and $\id N<\infty$, then $q_\fa(M,N)=q_\fa
(N)$ and $f_\fa(M,N)=f_\fa (\Hom_R(M,N))$.
\end{enumerate}
\end{corollary}

{\bf Proof.} i) The spectral sequence
$$H^{i}_{\fa}(\Ext_R^j(M,N))\underset{i}\Longrightarrow
H^{i+j}_{\fa}(M,N)$$  collapses at $j=p$, and so
$H^n_{\fa}(M,N)\cong H^{n-p}_{\fa} (\Ext_R^p(M,N))$ for all $n$.
Hence
$$
\begin{array}{ll} q_{\fa}(M,N)&=\sup\{n\in
\mathbb{N}_0:H^{n-p}_{\fa}(\Ext_R^p(M,N)) \text{ is  not Artinian }\}\\
&=\sup\{j+p:j\in\mathbb{N}_0  \text{ and } H^j_{\fa}(\Ext_R^p(M,N))
\text{ is not Artinian }\}
\\
&=q_\fa (\Ext^p_R(M,N))+p.
\end{array}$$
The similar argument shows that
$f_\fa(M,N)=f_{\fa}(\Ext^p_R(M,N))+p$.

ii) follows by i), Lemma 3.10 ii) and Lemma 3.10 v).

iii) By Lemma 3.10 vi), $\Gpd_NM=0$, and so in the light of the proof of
i), it turns out that $H^n_{\fa}(M,N)=H^n_{\fa}(\Hom_R(M,N))$ for
all $n$. Hence $f_\fa(M,N)=f_\fa (\Hom_R(M,N))$. On the other hand,
by \cite[Exersice 1.2.27]{BH}, $$\Ass_R(\Hom_R(M,N))= \Supp_RM\cap
\Ass_RN=\Ass_RN,$$ and so $\Supp_R(\Hom_R(M,N))= \Supp_RN.$ By Theorem
3.7 iv), this shows that
$$q_{\fa}(M,N)=q_{\fa}(\Hom_R(M,N))=q_\fa (N). \  \  \Box$$


\end{document}